\documentclass[11pt,reqno]{amsart}
\usepackage[margin=1in]{geometry}
\usepackage[T1]{fontenc}
\usepackage{lmodern,microtype,needspace}
\usepackage{amsmath,amssymb,amsthm,booktabs,xcolor}
\usepackage[colorlinks=true, pdfstartview=FitH, linkcolor=blue, citecolor=blue, urlcolor=blue]{hyperref}
\newcommand{\En}{\mathbb E^n}
\newcommand{\R}{\mathbb R}
\newcommand{\Z}{\mathbb Z}
\newcommand{\gap}{\operatorname{gap}}
\newtheorem{theorem}{Theorem}[section]
\newtheorem{lemma}[theorem]{Lemma}

\theoremstyle{definition}
\newtheorem{remark}[theorem]{Remark}
\title{Avoiding short progressions in Euclidean Ramsey theory II}
\author{Hei Chun Leung}
\address{Department of Mathematics, Hong Kong University of Science and Technology, Clear Water Bay, Hong Kong}
\email{hcleungat@connect.ust.hk}
\author{Chi Hoi Yip}
\address{Department of Mathematics, Hong Kong University of Science and Technology, Clear Water Bay, Hong Kong}
\email{machyip@ust.hk}

\begin{document}
\begin{abstract}
Recently, Currier, Moore, and Yip developed a general framework for constructing colorings avoiding prescribed monochromatic arithmetic progressions in Euclidean Ramsey theory. Writing $\ell_m$ for an $m$-term arithmetic progression with unit spacing, they proved $\En\not\to(\ell_3,\ell_{20})$ and $\En\not\to(\ell_3,\alpha\ell_{6889})$ for every $\alpha>0$. We refine and simplify their framework, proving $\En\not\to(\ell_3,\alpha\ell_{22})$ for every $\alpha>0$. We further establish $\En\not\to(\ell_3,\alpha\ell_{21})$ except possibly for four values of $\alpha^2$, and $\En\not\to(\ell_3,\alpha\ell_{17})$ for all but finitely many $\alpha>0$. Our approach also recovers several earlier results and avoids some technical arguments used in previous works.
\end{abstract}
\maketitle

\section{Introduction}

Let $\En$ denote $n$-dimensional Euclidean space. For $m\ge2$,
let $\ell_m$ denote an $m$-term arithmetic progression with
consecutive terms at distance one. For $\alpha>0$, write
$\alpha\ell_m$ for the corresponding configuration with spacing
$\alpha$. We refer to the square of the distance between consecutive
terms as the squared spacing. We write $\En\to(K_1,K_2)$ if every red--blue coloring of
$\En$ contains either a red congruent copy of $K_1$ or a blue
congruent copy of $K_2$. Otherwise, we write
$\En\not\to(K_1,K_2)$.

Euclidean Ramsey theory originated in the series of papers by
Erd\H{o}s, Graham, Montgomery, Rothschild, Spencer, and
Straus~\cite{EGMRSS1,EGMRSS2,EGMRSS3}; see also~\cite[Chapter~11]{Handbook}.
Positive results for arithmetic progressions include
$\mathbb E^2\to(\ell_2,\ell_5)$ due to
Tsaturian~\cite{Tsaturian} and
$\mathbb E^3\to(\ell_2,\ell_6)$ due to Arman and
Tsaturian~\cite{ArmanTsaturian}.  More recently, Currier, Moore, and
Yip~\cite{CMYPlane} established
$\mathbb E^2\to(\ell_3,\ell_3)$.

Spherical colorings, in which the color of a point depends only on
its distance from the origin, have long been used to obtain negative
results in every dimension. Erd\H{o}s et al.~\cite{EGMRSS1}
used such a coloring to prove $\En\not\to(\ell_6,\ell_6)$.
Conlon and Fox~\cite{ConlonFox} used a random periodic coloring
to prove $\En\not\to(\ell_2,\ell_{10^{5n}})$.
Currier, Mody, Xie, and Zhang~\cite{CMXZ} subsequently improved the
bound in arbitrary dimension, proving in particular that
$\En\not\to(\ell_2,\ell_m)$ for $m>(5+o(1))^n$ as
$n\to\infty$. They also established
$\mathbb E^2\not\to(\ell_2,\ell_{6330})$.
For red three-term progressions, Conlon and Wu~\cite{CW} combined
spherical colorings with probabilistic methods to establish
$\En\not\to(\ell_3,\ell_{10^{50}})$.
F\"uhrer and T\'oth~\cite{FT} subsequently gave a deterministic
spherical coloring showing
$\En\not\to(\ell_3,\ell_{1177})$.
They also treated progressions with unequal spacings, proving that
for each $\alpha>0$ there is an integer $m(\alpha)$ such that
$\En\not\to(\ell_3,\alpha\ell_{m(\alpha)})$,
and asked whether $m(\alpha)$ can be bounded independently of $\alpha$.

Recently, Currier, Moore, and Yip~\cite{CMY} answered this question, proving
\begin{equation}\label{eq:6889}
 \En\not\to(\ell_3,\alpha\ell_{6889}), \qquad \forall \alpha>0.    
\end{equation}
They also proved 
\begin{equation}\label{eq:345}
\En\not\to(\ell_3,\ell_{20}),
 \qquad \En\not\to(\ell_4,\ell_{14}), \qquad \En\not\to(\ell_5,\ell_8).
\end{equation}
Their framework colors each point $x\in\En$ according to the residue
of $\lfloor d\|x\|^2\rfloor$ modulo $p$. Along an arithmetic progression $\{x_k\}$, the values $d\|x_k\|^2$ form a
quadratic sequence, so the coloring conditions reduce to determining
which residues modulo $p$ occur among the integers
$\lfloor d\|x_k\|^2\rfloor$.
For arbitrary red residue sets, \cite[Proposition~2.4]{CMY}
gives a finite criterion for these questions.
For their uniform result, they additionally use the fact that a
sufficiently long quadratic sequence with integral or nearly
integral leading coefficient contains a translate of the set of quadratic residues or of the set
of quadratic nonresidues modulo a suitable prime $p$
\cite[Section~4]{CMY}. Consequently, their argument requires considering blue progressions of length on the order of $p^2$.

We refine this framework by studying the largest circular gaps among
values of quadratic sequences modulo an integer $p\ge2$, which need
not be prime. This leads to a unified criterion for excluding red
and blue progressions and allows us to treat arbitrary ratios between
their spacings. Our proof uses elementary properties of circular gaps
and exact finite computations, avoiding the more technical arguments
involving Diophantine approximation and quadratic sequences in \cite{CMY,FT}. The same criterion also recovers the three equal-spacing results in \eqref{eq:345}; see Remark~\ref{rem:CMY}.

Our first theorem significantly improves the $6889$ bound in \eqref{eq:6889} to $22$.
\begin{theorem}\label{thm:main}
For every $\alpha>0$ and positive integer $n$,
\[
 \En\not\to(\ell_3,\alpha\ell_{22}).
\]
\end{theorem}

For all but finitely many $\alpha>0$, the blue progression can be shortened to $17$ terms.
\begin{theorem}\label{thm:short}
For every positive integer $n$,
\[
 \En\not\to(\ell_3,\alpha\ell_{17})
\]
for all but finitely many $\alpha>0$. More precisely, this holds
whenever $\alpha^2$ is irrational, or $\alpha^2=a/b$ in lowest
terms with $a+b/50\ge37$, or
\[
 \frac7{152}\le\alpha^2\le\frac{179}{3600}.
\]
\end{theorem}

We also obtain a $21$-term bound with four possible exceptional
values of $\alpha^2$.
\begin{theorem}\label{thm:21partial}
For every positive integer $n$ and every $\alpha>0$ satisfying
\[
 \alpha^2\notin
 \left\{\frac14,\frac19,\frac2{15},\frac34\right\},
\]
we have
\[
 \En\not\to(\ell_3,\alpha\ell_{21}).
\]
\end{theorem}

Section~\ref{sec:prelim} presents our new spherical-coloring
framework and the exact criteria used to verify its hypotheses.
The three theorems are proved in Section~\ref{sec:results},
and the finite computations are described in
Section~\ref{sec:verification}.

\section{A general framework}\label{sec:prelim}

We first give a criterion for a spherical coloring to avoid red and
blue arithmetic progressions. We then show how to verify this criterion
using exact computations and how to choose the two spacings in a
prescribed ratio.

\subsection{The spherical-coloring criterion}

We begin with an elementary identity for the squared norms of points
in an arithmetic progression.

\begin{lemma}\label{lem:quadratic}
Let $M\ge2$. If $x_0,\ldots,x_{M-1}$ form an arithmetic progression
with consecutive points at distance $\sqrt U$, then there are real
numbers $\beta$ and $\gamma$ such that
\begin{equation}\label{eq:quadratic}
 \|x_k\|^2=Uk^2+\beta k+\gamma \qquad(0\le k<M).
\end{equation}
In particular, three consecutive points $x,y,z$ with spacing
$\sqrt R$ satisfy
\begin{equation}\label{eq:second}
 \|x\|^2-2\|y\|^2+\|z\|^2=2R.
\end{equation}
\end{lemma}

\begin{proof}
Let $v=x_1-x_0$ be the common difference vector. Then
$x_k=x_0+kv$ and $\|v\|^2=U$, so
\[
 \|x_k\|^2=\|x_0+kv\|^2
 =Uk^2+2k\langle x_0,v\rangle+\|x_0\|^2.
\]
This proves \eqref{eq:quadratic} with
$\beta=2\langle x_0,v\rangle$ and $\gamma=\|x_0\|^2$.
For three consecutive points, applying \eqref{eq:quadratic} with
$U=R$ gives
\[
 \|x\|^2-2\|y\|^2+\|z\|^2
 =\gamma-2(R+\beta+\gamma)+(4R+2\beta+\gamma)=2R.\qedhere
\]
\end{proof}

For $T>0$ and $I\subseteq\R$, we write $a\in I\pmod T$ if
$a=x+mT$ for some $x\in I$ and some integer $m$.
Fix integers $p\ge2$ and $1\le s\leq p-1$, and color $x\in\En$ red if
\begin{equation}\label{eq:color}
 \lfloor\|x\|^2\rfloor\in\{0,1,\ldots,s-1\}\pmod p,
\end{equation}
and blue otherwise. We do not assume that $p$ is prime.

For a nonempty finite set $A\subseteq\R/p\Z$, let
$\gap(A)$ denote the length of the largest gap between
consecutive elements of $A$ in cyclic order. If $A$
consists of a single point, we set $\gap(A)=p$. For $M\ge2$ and $u\in\R$, define
\begin{equation}\label{eq:G}
 G_{p,M}(u)=\max_{\beta\in\R}
 \gap\bigl\{
 -uk^2-\beta k\pmod p:0\le k<M
 \bigr\}.
\end{equation}
For fixed $u$, the expression being maximized is continuous and $p$-periodic in $\beta$, so the maximum is attained.

The following lemma gives sufficient conditions for the coloring \(\eqref{eq:color}\) to avoid red and blue arithmetic progressions of prescribed lengths and spacings.

\begin{lemma}\label{lem:conditions}
For the coloring \eqref{eq:color}, let $r,M\ge2$ and $R,U>0$.
\begin{enumerate}
\item If $G_{p,r}(R)\le p-s$, then the coloring contains no red
$r$-term progression with spacing $\sqrt R$.
\item If $G_{p,M}(U)\le s$, then the coloring contains no blue
$M$-term progression with spacing $\sqrt U$.
\end{enumerate}
For red three-term progressions with spacing $\sqrt R$, the simpler condition
\begin{equation}\label{eq:redcondition}
 R\in\left[s,\frac p2-s\right]\pmod{p/2}
\end{equation}
also suffices, provided that $4s<p$. 
\end{lemma}

\begin{proof}
Let $x_0,\ldots,x_{L-1}$ be an arithmetic progression
with spacing $\sqrt{u}$. By Lemma~\ref{lem:quadratic},
there exist $\beta,\gamma\in\R$ such that
\[
\|x_k\|^2=uk^2+\beta k+\gamma
\qquad(0\le k<L).
\]
Put $t_k=-uk^2-\beta k$. By the definition of the
coloring, the point $x_k$ is red precisely when
\[
\gamma\in[t_k,t_k+s)\pmod p,
\]
and blue precisely when
\[
\gamma\in[t_k+s,t_k+p)\pmod p.
\]

For fixed $\beta$, the first family of intervals covers
$\R/p\Z$ if and only if the largest circular gap between
the points $t_0,\ldots,t_{L-1}$ is at most $s$.
Similarly, the second family covers the circle if and
only if the largest circular gap between
$t_0+s,\ldots,t_{L-1}+s$ is at most $p-s$.
Since translation preserves circular gaps, these two
families have the same largest gap.

Consequently, if $G_{p,L}(u)\le s$, every such
progression contains a red point. Likewise, if
$G_{p,L}(u)\le p-s$, every such progression contains
a blue point. Taking $(L,u)=(M,U)$ and $(r,R)$,
respectively, proves the first two assertions.

For the last assertion, suppose that $x,y,z$ form a red three-term
progression with squared spacing $R$. Write
\[
 \|x\|^2=X+\varepsilon_x,\quad
 \|y\|^2=Y+\varepsilon_y,\quad
 \|z\|^2=Z+\varepsilon_z,
\]
where $X,Y,Z$ are integers and
$\varepsilon_x,\varepsilon_y,\varepsilon_z\in[0,1)$.
By \eqref{eq:second},
\begin{equation}\label{eq:floorerror}
 |(X-2Y+Z)-2R|<2.
\end{equation}
Since all three points are red, $X-2Y+Z$ belongs, modulo $p$, to
\[
 \{0,\ldots,s-1\}-2\{0,\ldots,s-1\}+\{0,\ldots,s-1\}
 =\{-2(s-1),\ldots,2(s-1)\}.
\]
Under \eqref{eq:redcondition}, we have $2R \in [2s,p-2s] \pmod p$, at circular distance at least $2$ from the above set, contradicting \eqref{eq:floorerror}.
\end{proof}

\subsection{Exact evaluation of the circular-gap function}

To apply Lemma~\ref{lem:conditions}, we need to bound $G_{p,M}$.
The next lemmas allow us to compute it exactly at rational
coefficients and certify its bounds on intervals of real coefficients.

\begin{lemma}\label{lem:period}
Let $p,M\ge2$ be integers. For all real numbers $u,v$,
\begin{equation}\label{eq:halfperiod}
 G_{p,M}(u+p/2)=G_{p,M}(u)
\end{equation}
and
\begin{equation}\label{eq:lipschitz}
 |G_{p,M}(u)-G_{p,M}(v)|
 \le\left\lfloor\frac{(M-1)^2}{4}\right\rfloor|u-v|.
\end{equation}
\end{lemma}

\begin{proof}
For every integer $k$,
\[
 -\left(u+\frac p2\right)k^2-
       \left(\beta-\frac p2\right)k
 =-uk^2-\beta k-\frac{k(k-1)}{2}p
 \equiv-uk^2-\beta k\pmod p,
\]
since $k(k-1)$ is even. Maximizing over $\beta$ gives
\eqref{eq:halfperiod}.

For the second assertion, write $v=u+\theta$ and compare the points $-uk^2-\beta k$ with
\[
-(u+\theta)k^2-(\beta-(M-1)\theta)k.
\]
Their difference is $\theta k(M-1-k)$, where
\[
0\le k(M-1-k)\le
C_M:=\left\lfloor\frac{(M-1)^2}{4}\right\rfloor.
\]
Since these differences lie in an interval of length
$C_M|\theta|$, the largest circular gap changes by at most
$C_M|\theta|$. Hence
\[
G_{p,M}(v)\le G_{p,M}(u)+C_M|v-u|.
\]
Interchanging $u$ and $v$ completes the proof.
\end{proof}

\begin{lemma}\label{lem:finitegap}
Let $u\in\mathbb Q$. Then $G_{p,M}(u)$ is rational, and the
maximum in \eqref{eq:G} is attained at one of the values
\begin{equation}\label{eq:beta-list}
 \beta=\frac{pm}{j-i}-u(i+j),\qquad
 0\le i<j<M,\quad 0\le m<j-i.
\end{equation}
\end{lemma}

\begin{proof}
Fix $u$. The cyclic order of the points in \eqref{eq:G} can
change only at values of $\beta$ satisfying
\[
 -ui^2-\beta i\equiv-uj^2-\beta j\pmod p
 \qquad(0\le i<j<M).
\]
Between successive such values, every circular gap is affine in
$\beta$, so their maximum is convex and attains its maximum
at an endpoint. Thus it suffices to consider the values of $\beta$
obtained by solving the displayed congruence. These are precisely
the values in \eqref{eq:beta-list}, where $0\le m<j-i$ suffices
because $\beta$ is considered modulo $p$. If $u$ is rational,
all the resulting gaps are rational.
\end{proof}

Thus $G_{p,M}(u)$ can be determined from the finitely many
values of $\beta$ in \eqref{eq:beta-list}. Together with
\eqref{eq:lipschitz}, this gives rigorous interval bounds using
only exact rational arithmetic.

\begin{remark}\label{rem:CMY}
The criterion of Lemma~\ref{lem:conditions} also recovers \eqref{eq:345} in which the red and blue
progressions have the same spacing.
Taking $R=U=d$, exact evaluation using Lemma~\ref{lem:finitegap}
gives
\[
\begin{array}{c|c|c|c}
(r,M)&(p,s,d)&G_{p,r}(d)&G_{p,M}(d)\\ \hline
(3,20)&(29,7,7)&22&67/10\\
(4,14)&(29,9,10)&20&62/7\\
(5,8)&(5,2,2)&3&2
\end{array}
\]
In each row, $G_{p,r}(d)\le p-s$ and $G_{p,M}(d)\le s$.
Rescaling gives the corresponding results of \eqref{eq:345}.
\end{remark}

\subsection{Choosing the spacings}

We now explain how to choose spacings in a prescribed ratio.
Put $T=p/2$. The conditions of Lemma~\ref{lem:conditions} are
periodic with period $T$, by Lemma~\ref{lem:period} and
\eqref{eq:redcondition}. Suppose that the coloring \eqref{eq:color} avoids red
$r$-term progressions whenever $R\in[A,B]\pmod T$, and blue
$M$-term progressions whenever $U\in[c,d]\pmod T$. For a prescribed $\alpha>0$, it therefore
suffices to find $R,U>0$ such that
\begin{equation}\label{eq:match-spacings}
 R\in[A,B]\pmod T,\qquad U\in[c,d]\pmod T,\qquad
 U=\alpha^2R.
\end{equation}
Rescaling by $1/\sqrt R$ then gives a coloring avoiding
$\ell_r$ in red and $\alpha\ell_M$ in blue.

When $\alpha^2$ is rational, the following lemma gives an exact
criterion for \eqref{eq:match-spacings}.

\begin{lemma}\label{lem:intersection}
Let $T>0$, and let $q=a/b>0$, where $a,b$ are coprime
positive integers. Suppose that $[A,B]$ and $[c,d]$ are closed
intervals, with $0<A\le B<T$. There exist $k\in\{0,\ldots,b-1\}$
and $r\in[A,B]$ such that
\[
 R=Tk+r,\qquad U=\frac abR,\qquad
 U\in[c,d]\pmod T
\]
if and only if
\begin{equation}\label{eq:intersection}
 \left\lceil\frac{bc-aB}{T}\right\rceil
 \le\left\lfloor\frac{bd-aA}{T}\right\rfloor.
\end{equation}
In particular, such $R,U$ exist whenever
\begin{equation}\label{eq:covering}
 a(B-A)+b(d-c)\ge T.
\end{equation}
\end{lemma}

\begin{proof}
The condition $U\in[c,d]\pmod T$ means that for some integer $m$,
\[
  c+Tm\le\frac ab(Tk+r)\le d+Tm.
\]
On multiplying by $b/T$ and putting $N=ak-bm$, this is equivalent to
\[
 \frac{bc-ar}{T}\le N\le\frac{bd-ar}{T}.
\]
As $r$ varies in $[A,B]$, the intervals above have union
\[
 \left[\frac{bc-aB}{T},\frac{bd-aA}{T}\right].
\]
Since $\gcd(a,b)=1$, every integer $N$ can be written as
$ak-bm$ with $0\le k<b$. Hence the required spacings exist
if and only if this last interval contains an integer, which is
precisely \eqref{eq:intersection}. Its length is
$[a(B-A)+b(d-c)]/T$, yielding \eqref{eq:covering}.
\end{proof}

For the red three-term condition \eqref{eq:redcondition}, we may take
$[A,B]=[s,T-s]$ when $4s<p$. In this case
\eqref{eq:covering} becomes
\[
 a(T-2s)+b(d-c)\ge T.
\]

The following lemma handles irrational squared spacing ratios
using the density of irrational rotations.

\begin{lemma}\label{lem:irrational}
Let $T>0$, $0<A<B<T$, and $u_0>0$. If $q>0$ is irrational,
then there exist arbitrarily large positive numbers $R,U$ such that
\[
 R\in(A,B)\pmod T,\qquad
 U\in u_0+T\Z,\qquad U=qR.
\]
\end{lemma}

\begin{proof}
For $m\ge0$, put $U_m=u_0+mT$ and $R_m=U_m/q$. Then
\[
 \frac{R_m}{T}=\frac{u_0}{qT}+\frac{m}{q}.
\]
Since $1/q$ is irrational, the fractional parts of $R_m/T$
are dense in $\R/\Z$. Every tail of this sequence is also dense,
so $R_m\in(A,B)\pmod T$ for arbitrarily large $m$.
The corresponding $R_m,U_m$ have the required properties.
\end{proof}

\section{Proof of the main results}\label{sec:results}

We first prove the $17$-term theorem, which also handles the cases
where $\alpha^2$ is irrational in the other two results.
For rational ratios, we use the interval-intersection criterion
and the exact computations described in Section~\ref{sec:verification}.

\subsection{The 17-term bound}\label{sec:short}

\begin{proof}[Proof of Theorem~\ref{thm:short}]
Use the coloring \eqref{eq:color} with $p=37$ and $s=9$.
By Lemma~\ref{lem:conditions}, the coloring contains no red three-term
progression with squared spacing $R\in[9,19/2]\pmod{37/2}$.
We first show that it contains no blue $17$-term progression
with squared spacing $U\in J\pmod{37/2}$, where
\begin{equation}\label{eq:17interval}
 J=\left[\frac7{16},\frac{179}{400}\right].
\end{equation}
Indeed, exact evaluation by Lemma~\ref{lem:finitegap} gives
\begin{equation}\label{eq:17values}
\begin{aligned}
 G_{37,17}\left(\frac{11}{25}\right)&=\frac{221}{25},\\
 G_{37,17}\left(\frac{111}{250}\right)&=\frac{222}{25},\\
 G_{37,17}\left(\frac{67}{150}\right)&=\frac{134}{15}.
\end{aligned}
\end{equation}
Since the Lipschitz constant in Lemma~\ref{lem:period}
is $64$ when $M=17$, we have $G_{37,17}(U)\le9$
whenever at least one of the following holds:
\[
\left|U-\frac{11}{25}\right|\le\frac1{400},
\qquad
\left|U-\frac{111}{250}\right|\le\frac3{1600},
\qquad
\left|U-\frac{67}{150}\right|\le\frac1{960}.
\]
The corresponding intervals overlap consecutively,
and their union contains $J$. This proves the claim by Lemma~\ref{lem:conditions} and \eqref{eq:halfperiod}.

Suppose first that $\alpha^2=a/b$, where $a,b$ are positive
coprime integers satisfying $a+b/50\ge37$.
Since $[9,19/2]$ and $J$ have lengths $1/2$ and $1/100$,
respectively, we have
\[
 a\left(\frac{19}2-9\right)
 +b\left(\frac{179}{400}-\frac7{16}\right)
 =\frac a2+\frac b{100}\ge\frac{37}2.
\]
Lemma~\ref{lem:intersection}, applied with $T=37/2$,
therefore gives positive squared spacings $R,U$ satisfying
$U=\alpha^2 R$, with $R\in[9,19/2]\pmod{37/2}$ and $U\in J\pmod{37/2}$.
Lemma~\ref{lem:conditions} and rescaling yield the conclusion.

The interval of values of $\alpha^2$ stated in the theorem
also follows directly. If \[
 \frac7{152}\le\alpha^2\le\frac{179}{3600},
\]
then $[9\alpha^2,\frac{19}{2}\alpha^2]$ intersects $J$.
We may therefore choose $R\in[9,19/2]$ such that
$U=\alpha^2R\in J$.

If $q=\alpha^2$ is irrational, apply Lemma~\ref{lem:irrational}
with $T=37/2$, $(A,B)=(9,19/2)$, and $u_0=11/25$.
The resulting squared spacings $R,U$ satisfy the red- and
blue-avoidance conditions established above. Rescaling completes
this case.

Any positive rational number $a/b$ in lowest terms not covered by the inequality $a+b/50\ge37$ satisfies
\[
 1\le a\le36,\qquad 1\le b<50(37-a),\qquad\gcd(a,b)=1.
\]
There are finitely many such values, proving the $17$-term assertion.
\end{proof}

\subsection{The uniform 22-term bound}\label{sec:proof}
\begin{proof}[Proof of Theorem~\ref{thm:main}]
Since $17<22$, Theorem~\ref{thm:short} gives the result whenever
$q=\alpha^2$ is irrational. We may therefore suppose that
$q=a/b$ is rational, with $a,b$ positive and coprime.
We use the three choices of $(p,s)$ listed below. For each,
the last column gives an interval $I$ such that the coloring~\eqref{eq:color}
avoids red three-term progressions whenever $R\in I\pmod{p/2}$:
\begin{center}
\begin{tabular}{ccc}
\toprule
$(p,s)$ & $T=p/2$ & Interval $I$ with $R\in I\pmod T$\\
\midrule
$(17,4)$ & $17/2$ & $[4,9/2]$\\
$(29,6)$ & $29/2$ & $[6,17/2]$\\
$(33,8)$ & $33/2$ & $[8,17/2]$\\
\bottomrule
\end{tabular}
\end{center}

For the coloring $(p,s)=(29,6)$, the exact calculation gives
\begin{equation}\label{eq:firstgap}
 G_{29,22}\left(\frac12\right)=\frac{17}{3}.
\end{equation}
By the Lipschitz estimate and periodicity in Lemma~\ref{lem:period},
every $U$ satisfying
\[
 U\in\left[\frac12-\frac1{330},\frac12+\frac1{330}\right]\pmod{29/2}
\]
satisfies $G_{29,22}(U)\le6$. With $[A,B]=[6,17/2]$
in Lemma~\ref{lem:intersection}, the sufficient condition
\eqref{eq:covering} becomes
\begin{equation}\label{eq:largenumerator}
 \frac{5a}{2}+\frac b{165}\ge\frac{29}{2}.
\end{equation}
In particular, this applies whenever $a\ge6$.

Next we use $(p,s)=(33,8)$. The exact interval verification
in Section~\ref{sec:verification} establishes
\begin{equation}\label{eq:biginterval}
 G_{33,22}(U)\le8
 \qquad\left(\frac9{25}\le U\le\frac{16}{25}\right).
\end{equation}
Taking $[A,B]=[8,17/2]$ in \eqref{eq:covering}, we obtain
\begin{equation}\label{eq:largedenominator}
 \frac a2+\frac{7b}{25}\ge\frac{33}{2},
 \qquad\text{equivalently}\qquad25a+14b\ge825.
\end{equation}
Combining \eqref{eq:largenumerator} and
\eqref{eq:largedenominator}, we are left with
\begin{equation}\label{eq:remaining}
 1\le a\le5,\qquad 25a+14b<825,\qquad
 b\ge1,\qquad\gcd(a,b)=1.
\end{equation}
There are $187$ such pairs.

The exact criterion \eqref{eq:intersection}, applied with
$[A,B]=[8,17/2]$ and $[c,d]=[9/25,16/25]$, holds for $104$ of
these pairs. For the remaining $83$ pairs, we use the values $u$
listed in Table~\ref{tab:values}. 

\begin{table}[ht]
\centering
\caption{Values used in the finite verification.}
\label{tab:values}
\begin{tabular}{ccc@{\qquad}ccc}
\toprule
$p$ & $u$ & $G_{p,22}(u)$ & $p$ & $u$ & $G_{p,22}(u)$\\
\midrule
$17$ & $1$ & $4$ & $29$ & $7$ & $6$\\
$17$ & $2$ & $4$ & $29$ & $14$ & $17/3$\\
$17$ & $1/2$ & $4$ & $29$ & $1/4$ & $23/4$\\
$17$ & $2/3$ & $133/39$ & $29$ & $157/11$ & $731/132$\\
$17$ & $3/5$ & $39/10$ & $33$ & $15/4$ & $15/2$\\
 & & & $33$ & $23/4$ & $8$\\
\bottomrule
\end{tabular}
\end{table}

The values of $G_{p,22}(u)$ in the table are computed using
Lemma~\ref{lem:finitegap}. In each row, $G_{p,22}(u)\le s$,
where $s=4,6,8$ for $p=17,29,33$, respectively. Put
\begin{equation}\label{eq:delta}
 \delta=\frac{s-G_{p,22}(u)}{110}.
\end{equation}
Then $G_{p,22}(U)\le s$ for every
$U\in[u-\delta,u+\delta]$.
Applying Lemma~\ref{lem:intersection} to these eleven intervals
handles all $83$ remaining pairs. The exact tests are given in
Section~\ref{sec:verification}.
This completes the proof.
\end{proof}

\subsection{The 21-term bound}\label{sec:m21}

\begin{proof}[Proof of Theorem~\ref{thm:21partial}]
By Theorem~\ref{thm:short}, the conclusion holds whenever
$\alpha^2$ is irrational. We may therefore assume that
$q=\alpha^2=a/b$, with $a,b$ coprime positive integers.

For $(p,s)=(33,8)$, Lemma~\ref{lem:finitegap} gives
\[
 G_{33,21}\left(\frac{3}{10}\right)=\frac{45}{7}.
\]
Since the Lipschitz constant in Lemma~\ref{lem:period} is
$100$, the interval
\begin{equation}\label{eq:m21interval}
 I=\left[\frac{199}{700},\frac{221}{700}\right]
\end{equation}
satisfies $G_{33,21}(U)\le8$ for every $U\in I$.
For the same coloring, \eqref{eq:redcondition} gives
$R\in[8,17/2]\pmod{33/2}$. This interval and $I$ have lengths
$1/2$ and $11/350$, respectively. Thus the sufficient condition
\eqref{eq:covering} holds whenever
\begin{equation}\label{eq:m21broad}
 \frac a2+\frac{11b}{350}\ge\frac{33}{2},
 \qquad\text{equivalently}\qquad 175a+11b\ge5775.
\end{equation}

There are exactly $5262$ pairs of coprime positive integers
$(a,b)$ with $175a+11b<5775$. The computation described in
Section~\ref{sec:verification} produces $728$ further intervals
on which $G_{p,21}(U)\le s$. Applying the exact criterion
\eqref{eq:intersection} to these intervals handles $5254$ of
the $5262$ pairs, leaving the eight values
\begin{equation}\label{eq:m21eight}
\left\{\frac14,\frac18,\frac19,\frac1{10},
        \frac2{11},\frac2{15},\frac2{23},\frac34\right\}.
\end{equation}

For four of these eight fractions, the following choices satisfy
$U=qR$, with $R\in[s,p/2-s]\pmod{p/2}$ and
$G_{p,21}(U)\le s$:
\begin{center}
\begin{tabular}{@{}cccccc@{}}
\toprule
$q$ & $p$ & $s$ & $R$ & $U$ & $G_{p,21}(U)$\\
\midrule
$2/23$ & $5$ & $1$ & $1$ & $2/23$ & $203/207$\\
$1/8$ & $9$ & $2$ & $23/2$ & $23/16$ & $2$\\
$1/10$ & $33$ & $8$ & $115/2$ & $23/4$ & $8$\\
$2/11$ & $5$ & $1$ & $27/2$ & $27/11$ & $75/77$\\
\bottomrule
\end{tabular}
\end{center}
The values in the last column are computed using
\eqref{eq:beta-list}. For each row, Lemma~\ref{lem:conditions}
shows that the coloring \eqref{eq:color}, with the listed
values of $p$ and $s$, contains neither a red three-term progression
with spacing $\sqrt R$ nor a blue $21$-term progression with
spacing $\sqrt U$. Rescaling by $1/\sqrt R$ completes the proof.
\end{proof}

The argument leaves open the four cases
\[
 \alpha^2\in\left\{\frac14,\frac19,\frac2{15},\frac34\right\}.
\]
It remains open whether the $21$-term conclusion holds for these
ratios. The finite verification tests only the specified family of
colorings and intervals, and does not establish that these four
ratios are exceptional.

\section{Verification of the finite computations}\label{sec:verification}

We describe the finite computations used in the proofs so that
the numerical claims can be checked independently.
By Lemma~\ref{lem:finitegap}, for rational $u$ it suffices to
consider the finitely many values of $\beta$ in
\eqref{eq:beta-list}. For each such $\beta$, we sort the points
$-uk^2-\beta k\pmod p$ in cyclic order and compute their largest
gap. The maximum over these values is $G_{p,M}(u)$. This also gives
the rational gap values stated in Section~\ref{sec:results}. All computations
and comparisons below use exact rational arithmetic. The source
code and supporting data are available at~\cite{GoodbillyCode}.

\subsection{The 22-term verification}
Every entry of Table~\ref{tab:values}, as well as
$G_{29,22}(1/2)=17/3$, is obtained by the preceding procedure.
We also need
\[
 G_{33,22}(U)\le8  \qquad\left(\frac9{25}\le U\le\frac{16}{25}\right).
\]
Starting from $[L,H]=[9/25,16/25]$, put $v=(L+H)/2$.
The Lipschitz estimate gives
\[
 G_{33,22}(U)\le G_{33,22}(v)+55(H-L)
 \qquad (U\in[L,H]).
\]
If the right-hand side is at most $8$, the inequality holds
throughout $[L,H]$; otherwise, we bisect the interval and apply the
same test to each half. The procedure terminates with $32$ subintervals and requires
$63$ exact evaluations of $G_{33,22}$. Since the Lipschitz bound
is applied on every subinterval, this establishes the inequality
for all real $U\in[9/25,16/25]$.

For the remaining $187$ rational pairs in \eqref{eq:remaining},
Lemma~\ref{lem:intersection} applied to the interval in
\eqref{eq:biginterval} yields
the exact test
\begin{equation}\label{eq:test-big}
\left\lceil\frac{18b-425a}{825}\right\rceil
\le\left\lfloor\frac{32b-400a}{825}\right\rfloor.
\end{equation}
This inequality holds for $104$ pairs. For each of the other $83$, take a row
$(p,u,G_{p,22}(u))$ from Table~\ref{tab:values}, put
$\delta=(s-G_{p,22}(u))/110$, with $s=4,6,8$ according to
$p=17,29,33$, and test
\begin{equation}\label{eq:test-small}
\left\lceil\frac{2b(u-\delta)-a(p-2s)}{p}\right\rceil
\le\left\lfloor\frac{2b(u+\delta)-2as}{p}\right\rfloor.
\end{equation}
For every remaining pair, the inequality holds for at least one of the eleven rows.
Thus these exact checks cover every pair in \eqref{eq:remaining}.

\subsection{The 21-term verification}
For each integer $9\le p\le65$, set
$s=\lfloor(p-1)/4\rfloor$ and evaluate $G_{p,21}(j/10)$ for
$0\le j<5p$. By \eqref{eq:halfperiod}, these values of $u$ represent all
multiples of $1/10$ modulo $p/2$. For each $u=j/10$ with
$g=G_{p,21}(u)\le s$, the Lipschitz
estimate gives $G_{p,21}(U)\le s$ whenever
\[
 U\in\left[u-\frac{s-g}{100},\,u+\frac{s-g}{100}\right]\pmod{p/2}.
\]
The computation yields $728$ distinct intervals, each associated
with a specified choice of $p$ and $s$ and interpreted modulo $p/2$.
There are
$5262$ pairs of coprime positive integers satisfying $175a+11b<5775$.
For each pair and each interval $[c,d]$, we apply
\eqref{eq:intersection} with $[A,B]=[s,p/2-s]$ and $T=p/2$.
For exactly $5254$ pairs, at least one test holds, leaving the eight
fractions in \eqref{eq:m21eight}. The four additional cases in Subsection~\ref{sec:m21}
are verified using \eqref{eq:beta-list}.
The files \path{verify_uniform22_with_rational17.cpp},
\path{verify_gaps.cpp}, \path{verify_coverage.py}, and
\path{verify_four.cpp}, together with the two exact gap tables
and instructions for reproducing the computations, are available
at~\cite{GoodbillyCode}.

\section*{Acknowledgments}
The authors used ChatGPT (OpenAI) to assist with the exposition
and the checking of calculations. The mathematical ideas and overall
approach are due to the authors, who take responsibility for all
mathematical statements and computational results.

\end{document}